\documentclass[12pt,thmsa]{article}
\usepackage{amsfonts}
\usepackage{amssymb}
\newtheorem{teo}{Theorem}[section]

\newtheorem{lema}[teo]{Lemma}
\newtheorem{cor}[teo]{Corollary}
\newtheorem{proposition}[teo]{Proposition}
\newtheorem{theorem}[teo]{Theorem}

\newtheorem{obs2}[teo]{Remark}

\newtheorem{tea}{Theorem}[subsection]

\newtheorem{no2}[teo]{Note}

\newtheorem{no3}[tea]{Note}

\newcommand{\Gal}{{\rm Gal}}
\newcommand{\Frob}{{\rm Frob }}
\newcommand{\trace}{{\rm trace}}
\newcommand{\mod}{{\rm mod}}
\newcommand{\lcm}{{\rm l.c.m.}}

\newcommand{\Q}{\mathbb{Q}}

\newcommand{\PSL}{{\rm PSL}}
\newcommand{\SL}{{\rm SL}}
\newcommand{\GO}{{\rm GO}}
\newcommand{\GU}{{\rm GU}}
\newcommand{\PGL}{{\rm PGL}}

\newcommand{\GL}{{\rm GL}}

\newcommand{\Image}{{\rm Image}}
\newcommand{\F}{{\mathbb{F}}}
\newcommand{\Z}{{\mathbb{Z}}}
\newcommand{\PSp}{{\rm PSp}}

\newcommand{\PSU}{{\rm PSU}}
\newcommand{\SU}{{\rm SU}}

\newcommand{\GSp}{{\rm GSp}}

\newcommand{\cond}{{\rm cond}}
\newcommand{\End}{{\rm End}}

\begin{document}
\title{{\bf Geometric families of $4$-dimensional Galois
representations with generically large images
}
}

\author{Luis Dieulefait
\and N\'{u}ria Vila
\\
Dept. d'\'{A}lgebra i Geometria, Universitat de Barcelona;\\
Gran Via de les Corts Catalanes 585;
08007 - Barcelona; Spain.\\
e-mail: ldieulefait@ub.edu; vila@mat.ub.es.
  }
\date{\empty}

\maketitle

\vskip 2cm
Research partially supported by MCYT grant BFM2003-01898 \\


Running head: Four-dimensional Galois representations\\

Author to send proofs:\\

N\'{u}ria Vila\\
Dept. d'\'{A}lgebra i Geometria, Universitat de Barcelona;\\
Gran Via de les Corts Catalanes 585;
08007 - Barcelona; Spain.\\
vila@mat.ub.es.

\newpage
\vskip 0.5cm

%

\begin{abstract}
We
study compatible families of four-dimensional Galois
representations constructed in the \'{e}tale cohomology of a smooth projective variety.
 We prove a theorem asserting that the  images will be
generically large if certain  conditions are satisfied.
 We only consider
representations with coefficients in
an imaginary quadratic field. We apply our result
 to an example constructed by Jasper Scholten (see [Sc]), obtaining
  a family of  linear groups
 and one of  unitary
groups as Galois groups over $\mathbb{Q}$.\\
\end{abstract}
%
\vskip 0.5cm
\section{Introduction}

In this article we study  geometric compatible
 families of four-dimensional Galois representations with coefficients in a quadratic
  imaginary field $K$
 and their images. We will assume that the representations have four different
  Hodge-Tate weights.\\
 After introducing some of the main tools,  
 we study
all possible types
 of resi\-dual image and show that,
 under certain conditions,
 the residual images, and the images themselves, of these geometric
representations are ``as large as possible" for almost every prime (i.e., for all
 but finitely many primes):
a  linear or a unitary group, depending on the decomposition
type of the prime in $K$. Then, we do explicit computations
 for one
example constructed by J. Scholten and
we bound the finite exceptional set, i.e., the set of primes where the image
is not ``as large as possible"  with an explicit
small density set of primes. Furthermore, we prove that this
example verifies all conditions of the result of the previous section
and, as a consequence, for this example the images are ``as large as possible"
for almost every prime.
We stress the consequences of these results for
Inverse Galois Theory.\\
Results of generically large images similar to the one obtained in this article
 were obtained by the authors in [D-V] for compatible families of
 three-dimensional Galois representations.\\
 The first large image result of this kind was obtained by Serre for the case of
 elliptic curves without Complex Multiplication (cf. [S1]). Serre also
 obtained a large image result for the four-dimensional Galois representations
 attached to abelian surfaces with trivial endomorphism ring (cf. [S2]), in his
 result he considers only symplectic representations (he can do this by reducing
  to the case of principally polarized abelian varieties), while we will
 work with four-dimensional representations which are not symplectic.
 The first author also proved a large image result for the case of
 four-dimensional Galois representations attached to genus two Siegel
 modular forms (cf. [D1]), but in that case he also restricts to
 symplectic representations. \\
 
\begin{flushleft}
\begin{flushleft}

\end{flushleft}
\end{flushleft}

\section{Geometric representations}
A compatible family of geometric four-dimensional $\lambda$-adic Galois
representations was constructed by Scholten (see [Sc]), we will recall this construction in section 7. In this section we will describe a fairly general class of compatible families of $4$-dimensional Galois representations containing the one constructed by Scholten, all the results in this article hold for any family in this class.\\
Let us start with a smooth projective variety $S$ defined over $\Q$ and we consider the action of $G_\Q$ on the groups $H^3(S)_\ell := H^3_{et}(S_{\bar{\Q}}, \Q_\ell)$.\\
 Let $N$ be the product of the primes of bad reduction of $S$. This gives a compatible family of Galois representations unramified at primes not dividing $N$. Let $K$ be an imaginary quadratic field, which for simplicity (and according to the construction in section 7) we will assume from now on to be $\Q(\zeta)$, where $\zeta$ is a cubic root of unity.\\

 Assume that the Galois representations on $H^3(S)_\ell \otimes_{\Q_\ell} K_\ell$ are all reducible, and they contain as  subrepresentations  a compatible family $\{ \sigma_{S, \lambda} \}$ of $\lambda$-adic $4$-dimensional Galois representations, with $\lambda$ prime in $K$ dividing $\ell$. For this family it also holds that the ramification set is (at most) the set of prime factors of $N$. We will assume  that the traces at Frobenius elements of these $4$-dimensional Galois representations, which lie all in $K$, generate $K$ (i.e., that they are not all rational integers). We will also assume that these representations have four different Hodge-Tate numbers: $\{0,1,2,3 \}$, and that the determinant of each $\sigma_{S, \lambda}$ is just $\chi^6$, where $\chi$ denotes the $\ell$-adic cyclotomic character (this last condition can always be satisfied after twisting by a sui\-ta\-ble Dirichlet character unramified outside $N$). \\

 Thus, the families we are considering can be shortly described as Galois representations associated to a pure motive defined over $\Q$ with coefficients in a quadratic field $K$ and Hodge-Tate numbers $\{0,1,2,3 \}$ (purity follows from Deligne's proof of the Weil's conjectures), with the simplifying assumptions that $K= \Q(\zeta)$ and the determinants are just $\chi^6$.\\
 
 From the compacity of $G_\Q$ and the continuity of the representations $\sigma_{S, \lambda}$, it follows that (after a suitable conjugation) we can assume that the images are contained in $\GL(4, \mathcal{O}_\lambda)$, where $\mathcal{O}$ denotes the ring of integers of $K$. This implies that we can consider the residual representations $\bar{\sigma}_{S, \lambda} $ with values in
 $\GL(4, \mathbb{F}_\lambda)$, obtained by composing $ \sigma_{S, \lambda} $ with the naive map ``reduction mo\-du\-lo $\lambda$".\\

Observe that from purity  it follows that if we denote by $a_p$ the trace
of the characteristic polynomial of $\sigma_{S, \lambda}( \Frob \; p)$ for any  prime $p$
 of good reduction (and $p \neq \ell$), it
holds: $|a_p| \leq 4 p^{3/2}$.\\
 
The characteristic polynomials of $\sigma_{S, \lambda} (\Frob \; p)$ ($p \nmid \ell N$) will be
denoted
$$ x^4 - a_p x^3 + b_p x^2 - p^3 \bar{a}_p x + p^6 \quad \quad (*)$$
The particular form of these polynomials follows, as is well known, from purity. Recall that the values $a_p$ generate $K$.\\
We will also assume (once again, this property holds in Scholten example) that all quadratic coefficients $b_p$ are rational integers.\\

Remark: We require that the traces at Frobenius elements generate $K$, an imaginary quadratic field, because in the opposite case, if all traces are real, the representations would be self-dual, thus for all primes $\lambda$ the images would not be ``as large as possible" (in fact, the images would be contained in an orthogonal group). Observe that even with the assumption that the traces generate $K$, the particular form of the characteristic polynomials in $(*)$ (together with Cebotarev
 density theorem) imply that for any prime $\lambda$ inert in $K$ the image of $\sigma_{S, \lambda}$ will be contained in the unitary group $\GU (4, \mathbb{Z}_\ell )$.\\

By conductor of the representations
$\sigma_{S, \lambda}$, $\cond(\sigma_{S, \lambda})$,
 we mean the prime-to-$\ell$ part of the Artin
conductor.
It is known that the conductor of the representations
$\sigma_{S, \lambda}$  is bounded independently of $\ell$, because the representations are constructed from a smooth projective variety.
This follows from
 Deligne's corollary to  the results of de Jong (cf. [Be], Proposition
6.3.2).
We will denote by
 $c:=\lcm\{\cond(\sigma_{S , \lambda})\}_{\lambda}$,
   this finite  uniform bound, and we will call it
    ``conductor of the family of representations".

\section{Classification of subgroups of $ \GL (4, q)$ }

The classification
 of the maximal subgroups $\mathcal{G}$, up to conjugation,
  of $\GL(4, p^t)$ can be derived from the results in Kleidman and Liebeck (cf. [K-L]),
   as done in [K] and in
  [N-P], pag. 592
  (see also 
  [A4-EMS], pags. 173 to 177).
   We will only consider the case $t=1 $ or $2$. In the following
   proposition, we are also using, for $t=2$, 
   the classification of maximal proper subgroups
   of $\GU (4,p^t )$ (also done in [K]).

   \begin{proposition}

   \label{teo:maximal}
   Let $\mathcal{G}$ be a subgroup of $\GL(4, p^t)$ with $t=
   1$ or $2$, $p >3$, and suppose that $\mathcal{G}$ does not contain
   $\SL (4,p^t )  $. Let $Z$ be the subgroup of scalar matrices contained in $\mathcal{G}$. Then at least one of the following
   holds:\\
   1) $\mathcal{G}$ is reducible;\\
   2) $\mathcal{G}$ is imprimitive, that is, conjugate either to (a) a
   subgroup of the group $\GL(1, p^t) \wr S_4$ of monomial matrices, or
   to (b) a subgroup of $\GL(2, p^t) \wr C_2$ (here $C_2$ denotes a group of order
   $2$);\\
   3) $\mathcal{G}$ is conjugate to a subgroup of the semidirect product of $\GL (2, p^{2t})$ by $C_2$ acting by field automorphisms;\\
   4) (only for $t=2$) $\mathcal{G}$ is conjugate to a subgroup of $
    \mathbb{F}_{p^2}^* \cdot \GL (4, p)$;\\
    5) $\mathcal{G}$ normalizes an extraspecial group of order $2^5$
    and $\mathcal{G}/Z$ is isomorphic to a subgroup of the split
    extension $2^4 \cdot S_5$;\\
    6) $\mathcal{G}$ is conjugate to a subgroup of $\GO^+ (4, p^t)$,
    $\GO^- (4, p^t)$ or $\GSp (4, p^t)$;\\
    7) $\mathcal{G}$ is absolutely irreducible and primitive, and
    $\mathcal{G}/Z$ is isomorphic to one of the following: $A_5$,
    $S_5$, $A_6$, $S_6$, $A_7$, $\PSL (2, 7)$ or $\PSp(4,3)$;\\
    8) (only for $t=2$) $\mathcal{G}$ is conjugate to a subgroup of
     $\GU (4, p)$.\\
    Furthermore, if $t=2$, $\mathcal{G}$ is contained in $\GU(4,p)$ (as in case (8))
     and $\mathcal{G}$ does not fall in any of cases (1) to (7), then
     $\SU(4,p) \subseteq \mathcal{G} \subseteq \GU (4,p)$.
\end{proposition}

\section[Etale and Crystalline cohomologies ]{Comparison
 of \'{E}tale and Crystalline cohomologies  and the
action of inertia}
As a consequence of the result of Fontaine-Messing in [F-M] giving a
(cano\-ni\-cal, functorial) isomorphism between the \'{e}tale cohomology and
the crystalline cohomology for a variety $S$ having good reduction at
$\ell$ and with $\dim (S)< \ell$, the following theorem, conjectured by
Serre in [S2], holds:\\
\begin{theorem}
\label{teo:Serrefm}
Let $V$ be a $ G_\mathbb{Q}$-stable lattice in
$ H_{et}^m (S_{\overline{\mathbb{Q}}
}, \mathbb{Q}_\ell)$ and $ \tilde{V}$ the
semi-simplification of $ V/\ell V$ with respect to the action of $ I =
I_\ell$. Then, if $ \dim (S) < \ell $ and $ S$ has good reduction at
$\ell$, the action of $I$ on each irreducible component $M$ of $\tilde{V}$ is given by a character
$ \Psi: I_t \rightarrow \mathbb{F}_{q}^*$ with $ q = \ell^n , n = \dim \; M $,
such that $$ \Psi = \Psi_f^{d_0 + \ell d_1 + ... + \ell^{n -1} d_{n-1}}$$ with $
\Psi_f$ fundamental character of level $n$ and the exponents verify:
$ 0 \leq d_i \leq m$. \\
\end{theorem}
Moreover, since we know that the Hodge-Tate numbers of our family of
four-dimensional representations are $\{ 0,1,2,3 \}$, if $\ell > 3$
 we can be more specific (cf. [F-L], theorem 5.3): the
digits $d_i$ appearing in the exponents above (for $M$ ranging over all irreducible
 components) must agree with the Hodge-Tate numbers. For example, if $\tilde{V}$
has, with respect to the action of $I_\ell$, an irreducible component $M_1$
of dimension one and another $M_2$ of dimension three, then the action
on $M_1$ is given by $\chi^a$, and the action on $M_2$ by
$\psi_3^{b + c \ell + d \ell^2}$, and it must hold: $\{ a,b,c,d\} = \{ 0, 1, 2,
3\}$. \\
Now that we know that the exponents of the action of $I_\ell$ through
fundamental characters are fixed, observe that there are many
possibilities for this action since it will depend on the dimension of the
irreducible components, and also on permutations of the four exponents.
Nevertheless, it is easy to see that in any case the  image
of inertia can not be too small. More precisely, an easy computation
gives the following:
\begin{lema}
\label{teo:inerciaislarge} Consider a four-dimensional crystalline
representation $\sigma_\lambda$ with Hodge-Tate numbers $\{ 0, 1, 2, 3\}$.
Then,
 if $\ell >3$, the
exponents of the fundamental characters giving the action of the (tame)
inertia subgroup $I_\ell$ of the residual representation $\bar{\sigma}_\lambda$
are also $\{ 0,1,2,3\}$. Let $P_I :=
\mathbb{P}(\bar{\sigma}_\lambda|_{I_\ell}^{s.s.})$
be the projectivization of the image of $I_\ell$. Then $P_I$ is a
cyclic group with order greater than or equal to $\ell - 1$.
\end{lema}

\section[Determination
 of the images Sa]{Determination
 of the images of the geome\-tric Galois representations }
We will consider a geometric compatible
family of four-dimensional Galois representations $\sigma_{S,\lambda}$ as the one described
in section 2.\\
We  will prove that the image of the residual
Galois representation
 is a  linear
group
 for almost every prime
decomposing in $\mathbb{Q}(\zeta)$, i.e., $ \ell \equiv 1 \pmod 3$,
and
a unitary group for  almost every inert prime $\ell \equiv 2 \pmod
3$,  whenever certain general conditions that
we will describe later are satisfied.

\subsection {Reducible representations}
\subsubsection{With a one dimensional component}

  Recall that $N$ denotes the product of the primes of bad reduction of
 $S$.
Suppose that $ \ell $ is a prime not in $S$, $\ell >3$, such that $
\bar{\sigma}_{S , \lambda}$
is reducible (not necessarily over the residue field $\F_\lambda$)
 with a one dimensional component.
In this case there is a character $ \mu$ unramified outside $ N \ell$
with image in $\bar \mathbb{F}^* _{\lambda}$ such that $ \mu (p)$ is a root
of $ x^4 - a_p x^3 + b_p x^2 - p^3 \bar{a}_p x + p^6$, for every $ p \nmid \ell N$.\\
Let $I := I_\ell$.
Using the description of $ \bar{\sigma}_{S , \lambda} | _I$ we know that
$ \mu = \chi^i \varepsilon$, with $\chi$ the $\mod \; \ell$ cyclotomic character,
 $ i = 0, 1, 2, 3$ and $ \varepsilon$ a
character unramified outside $N$. Clearly,
$\cond(\varepsilon)|\cond(\bar{\sigma}_{S , \lambda})|c$.\\
Using the description of $ \bar{\sigma}_{S , \lambda} | _I$ we can also
see that reducibility must necessa\-ri\-ly occur over $\F_\lambda$. If not,
there should be (at least) two one dimensional components, being each
conjugate of the other by some element of $\Gal( \bar{\F}_\lambda /
\F_\lambda)$. But on the other hand, lemma \ref{teo:inerciaislarge} implies that
the restriction of these two
one-dimensional components to $I$ must give two different powers
$\chi^i$ and $\chi^j$ of the cyclotomic character; if they where Galois
conjugated we would conclude that $i=j$ which is a contradiction.\\
Thus, the character $\epsilon$ has image contained in $\F^*_\lambda$.\\

At this point, we are led to introduce the following:\\

\bf{Condition 1)}: \it{Assume that there exists a prime $p$ not in $S$ such
that none of the roots of  $Pol_p(x)$, the characteristic polynomial of
${\sigma}_{S , \lambda} (\Frob \; p)$, is a number of the form $\eta
p^i$, where $\eta$ denotes an arbitrary root of unity and $i \in \{ 0,1,2,3
\}$}.\\

\rm
Under condition 1, we prove that there is a finite
number of  primes $\lambda$ such that the residual $\mod \; \lambda$
representation is reducible with a one dimensional component
$\mu_\lambda = \varepsilon_\lambda \chi^i$.\\
Assume on the contrary that this is the case for an infinite set of primes
 $\Lambda$. Since the characters $\varepsilon_\lambda$
  appearing in the one dimensional
 components of all these residual representations will have their
 conductors uniformly bounded by $c$, we can assume (by restricting to some
  infinite subset
  $\Lambda'$ of $\Lambda$) that there is a
  fixed character $\epsilon$ of conductor dividing $c$ with values in
  $\mathbb{C}^*$ such that the reduction of $\epsilon$ modulo $\lambda$
  gives $\varepsilon_\lambda$, for every $\lambda \in \Lambda'$. Of
  course, we can also assume (shrinking $\Lambda'$ if necessary) that
  the exponent of the cyclotomic character in $\mu_\lambda$ is
  independent of $\lambda \in \Lambda'$. \\
  Now take $p$ a prime as in condition 1), and if necessary delete $p$
  from $\Lambda'$. Then, for every $\lambda \in \Lambda'$ the polynomial
  $Pol_p(x)$ (which is
  independent of $\lambda$) when reduced modulo $\lambda$ has
  the reduction of $\epsilon(p) p^i$ modulo $\lambda$ as a root. \\
  Let $d$ be the order of $\epsilon(p)$. Then the above congruences
  implies that the resultant of the polynomials $Pol_p(x)$ and $x^d- p^{id}$
  is divisible by $\lambda$, for every $\lambda$ in $\Lambda'$: Thus,
  since this set of primes is infinite, we conclude that the resultant
  is $0$, and that the
  polynomial $Pol_p (x)$ has a root of the form: $ p^i \cdot \mbox{root of
  unity}$, contrary to condition 1). This proves the finitude of the
  primes in the resi\-dual redu\-cible case with one dimensional
  component, assuming condition 1).\\

\subsubsection{With two two-dimensional irreducible components}
This case can be treated in a similar way that the previous one, but
considering the exterior power $\rho_\lambda := \wedge^2 (\sigma_{S,
\lambda})$. The determinant of the
two-dimensional irreducible components of $\bar{\sigma}_{S,\lambda}$
 give rise to one dimensional components of $\bar{\rho}_\lambda$.\\
Using the information of the action of $I$ given in lemma
\ref{teo:inerciaislarge}, we see that a one dimensional component
 of $\bar{\rho}_\lambda$ must be of the form $\varepsilon \chi^i$, with
 the conductor of $\varepsilon$ dividing $c$ and the exponent $i \in \{
 1,2,3,4,5\}$. Again, using the restriction to $I$ one can easily see
 that the irreducible components can not be conjugated to each other, and so
  they must necessarily be defined over $\F_\lambda$. In particular,
  the image of $\varepsilon$ must be contained in $\F^*_\lambda$.\\
  Remark: Observe that, since we have shown that reducibility must necessa\-ri\-ly
  occur over the field of definition, this eliminates case 3) of the
  classification given in section 3.\\

  The natural condition to control this case is the following:\\

 \bf{ Condition 2)}: \it{Assume that there exists a prime $p$ not in $S$ such
that none of the roots of  $Q_p(x)$, the characteristic polynomial of
$\wedge^2({\sigma}_{S , \lambda}) (\Frob \; p)$, is a number of the form $\eta
p^i$, where $\eta$ denotes an arbitrary root of unity and $i \in
\{1,2,3,4,5
\}$. }\\

\rm
As before, if we assume that this case occurs for infinitely many primes, we get
 congruences modulo infinitely many primes between $Q_p (x)$ and $x^d -
p^{id}$, where $i \in \{1,2,3,4,5\}$ is fixed and $d$ is a fixed divisor of
$\varphi(c)$ (here $\varphi$ denotes Euler function), thus we conclude
that these two polynomials have a common root. For a prime $p$
verifying condition 2) this gives a contradiction.\\

\subsection{Imprimitive irreducible case}

\subsubsection{Monomial Case}

Consider first the monomial case: the image of $\bar{\sigma}_{S,\lambda}$ is a group $G$ having as normal subgroup the image $H$ of
 a one-dimensional representation, with quotient $G/H := U$ a subgroup of $S_4$.\\
Observe that since the largest cicle in $S_4$ has order $4$ and the projectivization of the image of the inertia subgroup $I$
is cyclic of order at least $\ell-1$ (by lemma \ref{teo:inerciaislarge}), if $\ell > 5$
the quotient $G/H$ gives an extension of $\mathbb{Q}$ with Galois group $U$
 unramified at $\ell$, thus ramifying only at primes in $N$.\\
 Remark: In [D-V] a similar situation was solved: the case of monomial
 three-dimensional Galois representations.\\

 We divide in two cases: \\
 a) $U$ contains odd permutations\\
 b) $U \subseteq A_4$\\

 a) In this case, $U$ has a quotient isomorphic to $C_2$, an element of $U$ corresponds to an odd permutation if and only if its image through this quotient is $-1$. This gives us a quadratic character $\psi$ corresponding to a quadratic extension of $\Q$ unramified outside $N$, which is a quotient of $\bar{\sigma}_{S, \lambda}$. An element $g$ of the Galois group of $\Q$ has  $\psi(g) = -1$ if
and only if the matrix  $\bar{\sigma}_{S, \lambda}(g)$ becomes diagonal after an odd permutation of its columns. It is easy to see that all such matrices have a property that can be read in their characteristic polynomials $Pol_g(x)$ (as a polynomial with coefficients in $\mathbb{F}_\lambda$): $Pol_g(x)$ has (at least) one pair of opposite roots $\alpha, -\alpha$. Taking $g = \Frob \;p$, $p \nmid N \ell$, we now translate this property in terms of the coefficients of $Pol_p(x)$:
$$ p^3 (a_p^2 + \bar{a}_p^2) \equiv a_p \bar{a}_p b_p \pmod{\lambda} \qquad \qquad (*)$$
(in general, for a polynomial $x^4+t_1 x^3 + t_2 x^2 + t_3 x + t_4$, the property of having two opposite roots gives the relation:
$t_3^2 + t_4 t_1^2 = t_1 t_2 t_3$).\\
Thus, we have seen that if the image of  $\bar{\sigma}_{S, \lambda}$
 falls in case a), and $\ell>5$, there should exist a quadratic
 character $\psi$ unramified outside $N$ such that for every prime
$p \nmid N \ell$ with $\psi(p) = -1$, the relation  (*) is satisfied. Therefore, using the principle that an equality of two algebraic integers is equivalent to infinitely many congruences:
$$  A=B \Longleftrightarrow A \equiv B \pmod{\lambda} , \mbox{for infinitely many prime ideals} \; \lambda $$
 we just have to impose the following condition to make sure that case a) can only occur for finitely many primes $\lambda$:\\

 \bf{Condition 3)}: \it{For every quadratic character $\psi$ unramified outside $N$ there exists a prime $p \nmid N$ with $\psi(p)=-1$ and
 $$ p^3 (a_p^2 + \bar{a}_p^2) \neq a_p \bar{a}_p b_p. $$}

\rm
b) In this case, we use the fact that $A_4$ has a normal subgroup $T$ isomorphic to $C_2 \times C_2$ with quotient
$A_4 / T \cong C_3 $. Since $U$ is a subgroup of $A_4$, if we denote $T' = T \cap U$, $T'$ is a normal subgroup of $U$
 with quotient $U / T' \cong C_3 $ or $T' = U$. We will assume that we are not in the case $T' = U$, because in this case
 the representation would either be reducible, or with image contained in
 $\GL(2, p^t) \wr C_2$, a case that we will consider later.\\
 Thus, we have a $C_3$ extension which is a quotient of our representation: this gives us a cubic character $\phi$ corresponding to a cubic extension of $\Q$ unramified outside $N$. An element $g$ of the Galois group of $\Q$ has $\phi(g) \neq 1$ if and only if the matrix  $\bar{\sigma}_{S, \lambda}(g)$ becomes diagonal after a permutation $P \in U$ of its columns such that the projection
 $\bar{P} \in U/T' \cong C_3$ is not the identity. It is easy to see that all such matrices have a property that can be read in
  their characteristic polynomials: they have three roots $\alpha, \beta, \gamma$ such that $\alpha + \beta + \gamma = 0$.\\
  Taking $ g = \Frob \; p $, this condition translates into the following restriction for the coefficients of $Pol_p(x)$:
  $$ a_p^2 b_p + p^6 \equiv p^3 \bar{a}_p a_p \pmod{\lambda} \qquad \qquad (**)$$
  (in general, for a polynomial $x^4+t_1 x^3 + t_2 x^2 + t_3 x + t_4$, the property of having three roots whose sum is $0$ gives the relation:
$t_1^2 t_2 + t_4 = t_3 t_1$).\\
Thus, if the image of the residual $\mod \; \lambda$ representation falls in case b) and $\ell >5$, there is a cubic character $\phi$ unramified outside $N$ such that for every prime $p \nmid N \ell$ with $\phi(p) \neq 1$, the relation (**) holds.\\
Obviously, the right condition that must be imposed to guarantee that this case can only occur for finitely many primes is the following:\\

\bf{ Condition 4)}: \it{For every cubic character $\phi$ unramified outside $N$ there exists a prime $p \nmid N$ with $\phi(p) \neq 1$ and  $$ a_p^2 b_p + p^6 \neq p^3 \bar{a}_p a_p .$$}
 \rm
 Observe that since the right hand side is a rational integer, it is enough for this to hold to have $a_p^2 b_p \not\in \mathbb{R}$.

 \subsubsection{Case 2)b) in the classification}
 This case already appeared in the study of symplectic representations done in [D1] and [D2], the required condition that the image of inertia must be not too small holds for every $\ell > 3$ by lemma \ref{teo:inerciaislarge}. Let us recall the result: if the image of  $\bar{\sigma}_{S, \lambda}$ falls in this case and $\ell >3$, then there is a quadratic character $\mu$ unramified
 outside $N$ such that for every prime $p \nmid N \ell$ with $\mu(p) = -1$, it holds:
 $$ a_p \equiv 0 \pmod \lambda.$$
 Observe that the situation is similar to the case of ``dihedral image" of two-dimensional Galois representations, studied by Serre and Ribet.\\
 We impose the following condition to ensure that this case does not occur for infinitely many primes:\\

 \bf{Condition 5)}: \it{For every quadratic character $\mu$ unramified outside $N$ there exists a prime $p \nmid N$ with $\mu(p)=-1$ and  $ a_p \neq 0.$}

\rm
 \subsection{Orthogonal and Symplectic groups}
If the image is contained in an orthogonal or symplectic group (case 6 of the classification) we use the following common property of these groups: the roots of the characteristic polynomial of any matrix come in reciprocal pairs. In particular, if we take
$Pol_p (x)$ the characteristic polynomial of $\bar{\sigma}_{S, \lambda} (\Frob \; p)$ for an unramified prime $p$, its roots will be
$\alpha, \beta, c / \alpha , c/ \beta$, with $c = \pm p^3$. This translates into the following condition on the coefficients of
$Pol_p(x)$: for every $p \nmid N \ell$
$$ \bar{a}_p \equiv \pm a_p \pmod{\lambda}. $$

 We impose the following condition to ensure that this case does not
  occur for infinitely many primes:\\

 \bf{Condition 6)}: \it{There exists a prime $p \nmid N$ such that: $a_p \neq \pm \bar{a}_p$.}\\

\rm
 Observe that this is equivalent to impose that $a_p^2 \not\in \mathbb{R}$.\\

\subsection{Case 4) in the classification}
We apply the same technique used to treat the corresponding case for symplectic
representations in [D1]. In this case, for any prime $p \nmid N \ell$, there should
exist a constant $c \in \F_{\ell^4}^*$ such that the reduction $\mod \; \lambda$ of
the polynomial:
$$ x^4 - c a_p x^3 + c^2 b_p x^2 + c^3 p^3 \bar{a}_p x + c^4 p^4 $$
a priori with coefficients in $\F_\lambda$,
 has coefficients in $\F_\ell$. This implies in parti\-cu\-lar that:
 $$ c^2 b_p \in \F_\ell ; \; \; c^2 a_p^2 \in \F_\ell $$
If we impose $\ell \nmid b_p$, (recall that, by assumption, $b_p \in \mathbb{Z}$),
we conclude first that $c^2 \in \F_\ell^*$, and then: $a_p^2 \in \F_\ell$.\\
Thus, any prime $\ell$ falling in this case must verify: $\ell \mid b_p$ or $a_p^2
\in \F_\ell$.\\
Clearly, the following condition implies that this case can only occur for finitely
many primes:\\

\bf{Condition 7)}: \it{There exists a prime $p \nmid N$ such that $b_p \neq 0$ and $a_p^2$
generates $\Q(\zeta)$.}\\

\rm
\subsection{Exceptional cases}
The case of ``exceptional groups", i.e., those listed as cases 5) and 7) in the
classification, can be discarded, for $\ell$ sufficiently large, just observing that
thanks to lemma \ref{teo:inerciaislarge} we know that the projective image of
inertia gives a subgroup of the image of order at least $\ell -1$ (see [D-V] and [D1]
for a similar argument). Computing the maximal order of an element in each of the
(projectivizations of the) exceptional groups, we easily see that for $\ell > 11$
these cases can not occur.

\subsection{Unitariness}
Recall that we are considering representations which are pure (i.e., all roots of
every characteristic polynomial of Frobenius have the same absolute value). This already reflected
in the coefficients of the characteristic polynomial, and it also implies that
for primes $\lambda$ inert in $\Q(\zeta)$ the image of the
representation is contained in the unitary group $\GU(4, \mathbb{Z}_\ell)$. Thus,
considering  residual representations, it is clear that their images, for every
prime
 inert in $\Q(\zeta)$, are contained in $\GU(4, \ell)$.\\

 \subsection{Conclusion}
Using the classification given in section 3, we conclude that any compatible family
of  representations of $G_\Q$ with the properties described in section 2 (geometric, pure,
with Hodge-Tate weights $\{0, 1, 2, 3\}$, coefficients in $\Q(\zeta)$, quadratic
coefficients in $\mathbb{Z}$ and determinant
 $\chi^6$) will have residual image ``as large as possible" for almost every prime,
  i.e., a general
linear group or a unitary group depending on the decomposition type of the prime in
$\Q(\zeta)$, as long as conditions 1) to 7) are satisfied. Using a lemma of Serre (cf. [S], [S2]),
this also gives ``maximality" of the image of the $\lambda$-adic representations for these primes.
For simplicity, we state the result in terms of the images of the projectivizations $P(\sigma_\lambda)$.

\begin{theorem}
\label{teo:imagenes} Let $ \{ \sigma_\lambda \}$ be a compatible family of geometric
pure $4$-dimensional Galois representations as those described in section 2, with
coefficients in $\Q(\zeta)$. Then, if the explicit conditions 1), 2),...7) given in
this section are satisfied, the image of $\sigma_\lambda$ is ``as large as possible"
for almost every prime, i.e., the image of its projectivization $P(\sigma_\lambda)$ satisfies:
$$ \Image(P(\sigma_\lambda)) \supseteq \PSL(4, \mathbb{Z}_\ell)$$ if $\ell$ splits in $\Q(\zeta)$, and
$$\Image(P(\sigma_\lambda)) = \PSU(4, \mathbb{Z}_\ell) $$ if $\ell$ is inert in $\Q(\zeta)$.
\end{theorem}

Remark: In the above result, in the case of splitting primes, let us give the exact value of the image of the projective residual representation $P(\bar{\sigma}_\lambda)$ according to the restriction imposed by $\det(\bar{\sigma}_\lambda) = \chi^6 $, a square in 
$\mathbb{F}_\ell$. This depends on the value of $\ell$ modulo $4$ (because this determines whether or not the order of
 $\mathbb{F}_\ell^*$ is divisible by $4$), so we easily see that for split primes with image ``as large as possible" we have:\\
 $$ \Image(P(\bar{\sigma}_\lambda)) = \PSL(4, \mathbb{F}_\ell)$$
 if $\ell \equiv 3 \pmod{4}$, and
 $$  \Image(P(\bar{\sigma}_\lambda)) = \{x \in \PGL(4, \mathbb{F}_\ell) \; : \; \det(x) \in (\mathbb{F}_\ell^*)^2 \} $$
  if $\ell \equiv 1 \pmod{4}$

\section{Explicit determination of the finite exceptional set}
Assuming that the value $c$ of the conductor of a geometric family is known, the
discussion of the previous section not only gives a set of conditions gua\-ran\-teeing
generically large image, but also an algorithm to explicitly bound the finite set of
exceptional primes (for a family verifying the 7 conditions). Let us write this
algorithm, which is just a recollection of the results obtained in the different subsections of
section 5. From now on we assume that  $\ell \nmid N$, $\ell >11$. At each step, only
 finitely many exceptional primes are detected: this follows from the analysis done in
  the previous section, under the assumption that the family of
 representations verify all conditions 1) to 7).\\
 
\begin{itemize}
    \item 
 Step 1)
 Take a few primes $\{ p_r \}$  with $Pol_{p_r}(x)$ as in condition 1).
 For each of
them compute the resultant $R_{r,i}$ of the pair of polynomials:
$$Pol_{p_r}(x) , x^d - {p_r}^{id}$$
with $i=0,1,2$ and $d$ equal to the order of $p_r$ in $(\Z / c \Z)^*$. If $\ell$
 does not divide the resultants $R_{r,0}, R_{r,1}$ and $R_{r,2}$ for one of these
 primes $p_r$ (and $\ell \neq p_r$) then $\ell$ is not exceptional for step 1). \\
 
 \item Step 2) Take a few primes $\{ q_r \}$  with $Q_{q_r}(x)$ as in condition 2).
 For each of
them compute the resultant $R_{r,i}$ of the pair of polynomials:
$$Q_{q_r}(x) , x^d - {q_r}^{id}$$
with $i=1,2,3,4,5$ and $d$ equal to the order of $q_r$ in $(\Z / c \Z)^*$. If $\ell$
 does not divide the resultants $R_{r,1}, R_{r,2},...R_{r,5}$ for one of these
 primes $q_r$ (and $\ell \neq q_r$) then $\ell$ is not exceptional for step 2). \\

\item Step 3) For every quadratic character $\psi$ unramified outside $N$ take a few
primes $\{ p_r\}$ such that $\psi(p_r) = -1$ for all of them and they also satisfy
the rest of condition 3). Then, if for some of these $p_r$, the prime $\ell$ does
not divide:
$$ p_r^3 (a_{p_r}^2 + \bar{a}_{p_r}^2) - a_{p_r} \bar{a}_{p_r} b_{p_r} $$
and $\ell \neq p_r$, the prime $\ell$ is not exceptional ``with respect to $\psi$"
for step 3). Repeating the process for the finitely many possible $\psi$, we compute
all exceptional primes for step 3).
\\

\item Step 4) For every cubic character $\phi$ unramified outside $N$ take a few primes
$\{ p_r\}$ such that $\phi(p_r) \neq 1$ for all of them and they also satisfy the
rest of condition 4). Then, if for some of these $p_r$, the prime $\ell$ does not
divide:
 $$ a_{p_r}^2 b_{p_r} + p_r^6 - \bar{a}_{p_r} a_{p_r} $$
and $\ell \neq p_r$, the prime $\ell$ is not exceptional ``with respect to $\phi$"
for step 4). Repeating the process for the finitely many possible $\phi$, we compute
all exceptional primes for step 4). \\

\item Step 5) For every quadratic character $\psi$ unramified outside $N$ take a few
primes $\{ p_r\}$ such that $\psi(p_r) = -1$ for all of them and they also satisfy
the rest of condition 5) (i.e., $p_r \nmid N$ and $a_{p_r} \neq 0$). Then, if for
some of these $p_r$, the prime $\ell$ does not divide $a_{p_r}$ and $\ell \neq p_r$,
the prime $\ell$ is not exceptional ``with respect to $\psi$" for step 5). Repeating
the process for the finitely many possible $\psi$, we compute all exceptional primes
for step 5).\\

\item Step 6) Take a few primes $\{ p_r \}$  with $a_{p_r}$ as in condition 6), i.e., $p_r
\nmid N$, $a_{p_r}^2 \not\in \mathbb{R}$. Then, if for some of these $p_r$, the
prime $\ell$ does not divide $a_{p_r}^2 - \bar{a}_{p_r}^2$ and $\ell \neq p_r$, the
prime $\ell$ is not exceptional  for step 6).\\

\item Step 7) Take a few primes $\{ p_r \}$ verifying condition 7). Then, if for some of
these $p_r$ the prime $\ell$ verifies:
 $\ell \nmid b_{p_r}$, $a_{p_r}^2
\not\in \F_\ell$ and $\ell \neq p_r$, the prime $\ell$ is not exceptional for step
7). \\
\end{itemize}

Thus, 
we bound the set of exceptional primes with the finite set of primes which are
exceptional at some step of the algorithm, or divide $N$, or are smaller than
$13$.\\

Remark:
One can expect that by  sufficiently enlarging  the set of ``auxiliary primes" $\{
p_r \}$ at each step, one can get a bounding set which contains only the true
exceptional primes, and the primes dividing $N$ or smaller than $13$. The problem is that the size of such a ``sufficiently large"
 auxiliary set  may be (according to the existing effective versions of Cebotarev density theorem) too large to be suitable for computations.

\section{Explicit computations on the example}
We will now apply the former results to the non-selfdual 4-dimensional $\ell$-adic Galois representation given by Scholten (see [Sc]): we check that the image is ``as large as possible" for almost every prime, but since the value of the conductor of the family is unknown, we can not give a finite bound
 for the set of exceptional primes. Instead, we can bound the set of exceptional primes by a small density set of primes by a method already applied in [D-V].

\subsection{The example}
First, we summarize the construction of a compatible family of geometric four dimensional $\lambda$-adic Galois representation given by Scholten (see [Sc]). This family belongs to the class of families discussed in section 2 (in fact, it inspired the definitions therein).\\ 
The example is constructed using a subspace of
$H^3_{et}( S_{\overline{ \mathbb{Q}}}, \mathbb{Q}_\ell)$, where $S$ is  a  threefold obtained from a base change fibre product of two explicitly given elliptic surfaces. 

Let $X$ and $X'$ be the
  elliptic surfaces given by the Weierstrass equations:

$$X:\quad y^2=x^3-27(t^2+t+1)^2 x+18(3t^2+3t-1)(t^2+t+1)^2,$$

$$X':\quad y^2=x^3-3(3t+2) (243t^3+486t^2+324t+80)x-39366t^6-157464t^5$$
$$ -262440t^4-235224t^3-120528t^2-33696t-4048.$$

Let $\mathcal{E}$ and $\mathcal{E}'$ be the Kodaira-N\'{e}ron models of the cubic base change of the elliptic surfaces $X\rightarrow\mathbb{P}^1_t$ and $X'\rightarrow\mathbb{P}^1_t$ over the Galois cover $\varphi:\mathbb{P}^1_s\rightarrow\mathbb{P}^1_t$ given by 
$$t=\varphi(s)=\frac{s^3-3s^2+1}{3s^2-3s}.$$
The Galois group of $\varphi$ has order 3 and it is generated by $\phi:s\mapsto(s-1)/s$.

Let $S$ denote  the blowup at the singular points on $\mathcal{E}\times_{\mathbb{P}^1}\mathcal{E}'$.
From the construction, the threefold $S$ has bad reduction at the primes dividing $N=6$ and $\phi$ induces an automorphism on $S$ defined over $\mathbb{Q}$.
Let $\zeta$ denote a fixed primitive third root of unity. The automorphism $\phi$ induces a $G_\mathbb{Q}$-invariant linear automorphism on the $\ell$-adic cohomology, also denoted by $\phi$. Since $\phi$ has order 3, the $\ell$-adic $G_\mathbb{Q}$-module 
$H^3(S)_\ell:=H^3_{et}( S_{\overline{ \mathbb{Q}}}, \mathbb{Q}_\ell)$ decomposes as a sum of three eigenspaces over $\mathbb{Q}_\ell(\zeta)$. Let $W_\lambda$ denote the $\zeta$-eigenspace of the automorphism  $\phi$ restricted to $H^3(S)_\ell \otimes_{\mathbb{Q}_\ell}\mathbb{Q}_\ell(\zeta)$.
 In  [Sc] it is shown, assuming a statement about the generators of the $H^2_{et}( S_{\overline{ \mathbb{Q}}}, \mathbb{Q}_\ell)$,  that  the $\zeta$-eigenspace $W_\lambda$ has dimension 4 and that $$\dim H^{3,0}(W_\lambda)=\dim H^{2,1}(W_\lambda)=\dim H^{1,2}(W_\lambda)=\dim H^{0,3}(W_\lambda)=1.$$
So, we have a compatible family of $\lambda$-adic four dimensional
Galois representations $\sigma'_\lambda$ on $W_\lambda$,
a $\mathbb{Q}(\zeta)_\lambda$ vector space, with four different
Hodge-Tate numbers: $\{ 0,1,2,3 \}$. In general
$\sigma'_\lambda$ has to be twisted by a Dirichlet character $\epsilon$
to obtain $\epsilon \otimes \sigma'_\lambda = \sigma_\lambda$ with
$\det (\sigma_\lambda) = \chi^6$.\\
This is the compatible family of geometric Galois  representations we will consider:
$$ \sigma_{S,\lambda} : G_{\mathbb{Q}} \rightarrow \GL(W_\lambda).$$

On the other hand, we put again $a_p := \trace (\sigma_{S , \lambda}(\Frob \; p)) \in
 \mathbb{Q}(\zeta)$
for every $p \nmid \ell 6$, we have  
$|a_p| \leq 4 p^{3/2}$.
These traces can be determined using Lefschetz trace formula and counting points of $S$
over finite fields. 
The characteristic polynomials of $\sigma_{S ,\lambda} (\Frob \; p)$ will be
$$ x^4 - a_p x^3 + b_p x^2 - p^3 \bar{a}_p x + p^6 .$$
For all primes $p \nmid \ell 6$, the quadratic coefficients $b_p$ are rational integers.\\
 According with [Sc] we list the first characteristic polynomials: \\

$  5:  x^4+(\overline{13+10\zeta})x^3-5x^2+5^3(13+10\zeta)x+5^6 $ \\
$  7:  x^4+(\overline{7+3\zeta})x^3-189x^2+7^3(7+3\zeta)x+7^6 $ \\
$  11:  x^4+(\overline{21+2\zeta})x^3+517x^2+11^3(21+2\zeta)x+11^6 $ \\
$  13:  x^4+(\overline{70+77\zeta})x^3-1742x^2+13^3(70+77\zeta)x+13^6 $ \\
$  17:  x^4+(\overline{87+63\zeta})x^3-1802x^2+17^3(87+63\zeta)x+17^6 $ \\
$  19:  x^4-(\overline{8+81\zeta})x^3-4275x^2-19^3(8+81\zeta)x+19^6 $ \\
$  23:  x^4+(\overline{129+33\zeta})x^3+14536x^2+23^3(129+33\zeta)x+23^6 $ \\
$  29:  x^4+(\overline{186+86\zeta})x^3+16936x^2+29^3(186+86\zeta)x+29^6 $ \\

Remark: In this example, we have twisted by a cubic character
$$ \epsilon: G_{\mathbb{Q}} \rightarrow C_3 $$
of conductor $9$
to obtain $\det ( \sigma_{S , \lambda} ) = \chi^6$, $\chi$ the $\ell$-adic
cyclotomic character.

\subsection{Determination of the images of $\sigma_{S , \lambda}$}
Using the $8$ characteristic polynomials listed above it is easy to check that the family of Galois representations $ \sigma_{S,\lambda }$ satisfies the  conditions of theorem \ref{teo:imagenes}.
One easily checks that the characteristic polynomial at $5$ satisfies conditions 1) and 2) (all roots and products of two roots
 of this characteristic polynomial generate non-abelian extensions of $\Q$).
We consider all quadratic fields  unramified outside $6$ and we check that conditions 3) and 5) are satisfied for some prime $p$ inert in each of them. The same also holds for cubic characters and condition 4). Conditions 6) and 7) are easily checked . Thus, we conclude:

\begin{theorem}
\label{teo:imagenesejemplo} The images of the  
 $4$-dimensional Galois representations $ \sigma_{S,\lambda }$  are ``as large as possible"
for almost every prime, i.e., for almost every prime $\lambda$:
$$ \Image(P(\sigma_{S , \lambda})) \supseteq \PSL(4, \Z_\ell)$$ if $\ell$ splits in $\Q(\zeta)$ ($\ell \equiv 1 \pmod{3}$), and
$$\Image(P(\sigma_{S , \lambda})) = \PSU(4, \Z_\ell)$$  if $\ell$ is inert in $\Q(\zeta)$ ($\ell \equiv 2 \pmod{3}$).
\end{theorem}

We can not bound the set of exceptional primes with an explicit finite set (using the algorithm in section 6) because the value of the conductor of the family is unknown. In order to obtain a partial similar result, we will proceed as in [D-V] and take a fake value for this conductor: $c = 27 \cdot 64$. Recall that the value of $c$ is only needed in order to deal with the reducible cases in steps 1) and 2) of the algorithm (section 6), and as a matter of fact (see section 5.1) it is only needed in order to bound the conductors of the Dirichlet characters appearing in the residually reducible cases, which are all characters unramified outside $6$ and valued in 
$\F_\lambda^*$. \\
Therefore, what we really need is a bound for the conductors of these Dirichlet characters, and using the fact that they take values in $\F_\lambda^*$ we can bound these conductors after excluding a suitable set of primes (the density of this set tends to $0$ when the value for the fake conductor tends to infinity): this follows easily from the fact that $\F_\lambda^*$ is cyclic and has order
$\ell -1$ if $\ell$ is split in $\Q(\zeta)$ and $\ell^2 - 1$ otherwise. In particular, we easily check that the conductors of the Dirichlet characters are bounded by $27 \cdot 64$ if we exclude those primes in the following sets:\\
1) Splitting primes, $\ell \equiv 1 \pmod{3}$: $$D_1 := \{\ell \; : \; \ell \equiv 1 \pmod{27} \} \cup
\{\ell \; : \; \ell \equiv 1 \pmod{32} \} $$
2) Inert primes, $\ell \equiv 2 \pmod{3}$ : $$D_2 := \{\ell \; : \; \ell \equiv -1 \pmod{27} \} \cup
\{\ell \; : \; \ell \equiv \pm 1 \pmod{16} \} $$
Thus, from now on we only consider primes $\lambda$ such that $\ell$ is not in $D_1 \cup D_2 $, and this implies that we can apply the algorithm in section 6 taking $c = 27 \cdot 64$.\\

Remark: In terms of Dirichlet density, we are excluding $1/6$ of the splitting primes and $1/3$ of the inert primes.\\  

We executed the algorithm in Pari GP, using only the $8$ characteristic polynomials listed in the previous subsection, and we found no exceptional prime greater than $11$. 

 \begin{theorem}
\label{teo:imagenesejemplo2} The images of the  
 $4$-dimensional Galois representations $ \sigma_{S,\lambda }$  are ``as large as possible" (as described in the previous theorem)
 for every prime $\lambda$ such that $\ell >11$ and $\ell \not\in D_1 \cup D_2 $.
\end{theorem}

If we consider the projectivizations of the residual representations, for those primes $\lambda$ where the image is ``as large as possible" (and, for the case of splitting primes, we also impose $\ell \equiv 3 \pmod{4}$, according to the remark after theorem
 \ref{teo:imagenes}), we obtain as a corollary two families of classical groups over finite fields realized as Galois groups over $\Q$ . These groups were not previously known to be Galois groups over $\Q$ (cf. [M-M], [V]).

\begin{cor}
\label{teo:GaloverQ} Let $\ell$ be a prime greater than $11$ not in $D_1 \cup D_2$. Then, the following groups are Galois groups of a finite extension of $\Q$ unramified outside $6 \ell$:\\
1) If $\ell \equiv 7 \pmod{12}$: $\PSL(4, \ell)$.\\
2) If $\ell \equiv 2 \pmod{3}$: $\PSU(4, \ell)$.
\end{cor}

%


\begin{thebibliography}{}
%
%
\bibitem[A4-EMS]{RefJ}
Kostrikin, A.I., Shafarevich, I. R. (Eds), Algebra IV, Encyclopedia of Mathematical Sciences, vol. 37, Springer-Verlag, 1989.
\bibitem[Be]{RefJ}
 Berthelot, P., Alt\'{e}rations de vari\'{e}t\'{e}s alg\'{e}briques [d'apr\`{e}s A. J. de
 Jong], in ``S\'{e}minaire Bourbaki", vol. 1995/96, exp. 815, 273-311; Ast\'{e}risque, 1997.
\bibitem[D1]{RefJ}
 Dieulefait, L.,  On the images of the Galois representations
 attached to genus $2$ Siegel modular forms,  J. Reine Angew. Math. {\bf 553} (2002), 183-200.   
 \bibitem[D2]{RefJ}
  \_\_\_\_\_\_ , Explicit determination of the images of the
 Galois representations attached to abelian surfaces with $\End(A) =
  \mathbb{Z}$, Experiment. Math. {\bf 11} (2002), 503-512.   \bibitem[D-V]{RefJ}                                                                                  Dieulefait, L., Vila, N., On the images of modular and geometric three-dimensional 
      Galois representations, Amer. J. of Math. {\bf 126} (2004), 335-361.  
    \bibitem[F-L]{RefJ}
 Fontaine, J.M., Laffaille, G., Construction de repr\'{e}sentations
$p$-adiques, Ann. Scient. \'{E}c. Norm. Sup., $4^e$ s\'{e}rie, t. {\bf 15} (1982), 547-608.
\bibitem[F-M]{RefJ}
 Fontaine, J.M, Messing, W., $p$-adic periods and $p$-adic
etale cohomology, in ``Currents Trends in Arithmetical Algebraic
Geometry (Arcata, Calif., 1985)"; Contemporary Mathematics {\bf 67}
(1987), 179-207. 
\bibitem[K]{RefJ}
Kleidman, P. B., The subgroup structure of some finite simple groups, Ph. D. Thesis, 
Cambridge, 1986. 
\bibitem[K-L]{RefJ}
 Kleidman, P., Liebeck, M., \textit{The subgroup structure of the finite classical groups}, London Math. Soc. LNS 129,   Combridge University Press, 1990.
\bibitem[M-M]{RefJ}
 Malle, G., Matzat, B. H., \textit{Inverse Galois Theory}, Springer-Verlag, 1999.
 \bibitem[N-P]{RefJ}
 Neumann, P. M., Praeger, C. E., A recognition algorithm for special linear groups, Proc. London Math. Soc. (3) {\bf 65} (1992), 555-603.
\bibitem[Sc]{RefJ}
 Scholten, J., Mordell-Weil groups of elliptic surfaces and
Galois representations, Ph. D. Thesis, Rijksuniversiteit Groningen, 2000.
\bibitem[S]{RefB}
 Serre, J-P., \textit{Abelian $\ell$-adic representations and elliptic
curves}, Benjamin, 1968.
\bibitem[S1]{RefJ}
 \_\_\_\_\_\_ , Propri\'{e}t\'{e}s galoisiennes des points d'ordre
fini des courbes elliptiques, Invent. Math. {\bf 15} (1972), 259-331.
\bibitem[S2]{RefJ}
 \_\_\_\_\_\_ , Oeuvres, vol. 4, 1-55, Springer-Verlag, 2000.
  \bibitem[V]{RefJ}  
V\"{o}lklein, H., Braid group action via GL$_n(q)$ and U$_n(q)$ and Galois realizations, Israel J. Math. {\bf 82} (1993), 405-427. 
\end{thebibliography}
\end{document}